\documentstyle[12pt]{article} 
\begin{document} \begin{center}{\large {\bf Deformations of Semi-Euler
Characteristics}}\bigskip\\ {\sc George R. Kempf\\ The Johns Hopkins University
and\\ University of California, Riverside}\bigskip \bigskip \\ \end{center}

Let $f:X\rightarrow S$ be a proper smooth morphism of pure relative dimension
$n$ with connected fibers between Noetherian schemes.  Let $\xi$ be a locally
free coherent sheaf on $X$.  If $s$ is a point of $S$, we have the sheaf $\xi_s
= \xi|_{X_{s}}$ on the fiber $X_s$ of $f$ over $s$.\\

If $G$ is a coherent sheaf on a proper variety $Y$, the semi-Euler
characteristic
$\psi(G) = \displaystyle{\sum_{i \mbox{ {\scriptsize even}}}} \dim H^i(Y,G)$.\\

If $n$ is odd say $1+2m$, we will assume that we are given a non-degenerate
pairing $B:\xi \otimes \xi \rightarrow \omega_{X/S}$ such that $B$ is symmetric
if $n \equiv  1(4)$ or skew-symmetric if $n\equiv 3(4)$.  In this situation we
have\\

\noindent {\bf Theorem 1.}  The parity of $\psi (\xi_s)$ is locally constant on
$S$ if 2 is a unit in ${\cal O}_S$.\\

In characteristic zero this result appears in [1].  If $n=1$, the result
appears
in [2].  My proof uses that the rank of a skew-symmetric matrix is even and it
yields deeper results on the variation of $\dim H^i(X_s, \xi _s)$.\\

\noindent {\bf \S1.  The statement of the main results.}\\

A Grothendieck complex for $\xi$ is a complex $K^*:0\rightarrow K^0
\stackrel{\alpha^0}{\rightarrow }K^1 \stackrel{\alpha^1}{\rightarrow }\cdots
K^n
\rightarrow 0$ of free coherent sheaves of $S$ such that the $i$ cohomology of
$K^*|_s$ is isomorphic to $H^i(X_s, \xi_s)$ for all points $s$ in $S$.
Grothendieck has shown that such complexes always exist locally on $S$.\\

We say that the complex $K^*$ is special if it has the form \[0 \rightarrow K^0
\stackrel{\alpha^0}{\rightarrow }K^1 \rightarrow \cdots \rightarrow
K_m\stackrel{\beta}{\rightarrow }\check K_m \stackrel{(-1)^*\check
\alpha_{m-1}}{\rightarrow } \check K_{m-1}\rightarrow \cdots\rightarrow \check
K^1 \stackrel{(-1)^*\check \alpha^0}{\rightarrow }\check K^0\rightarrow 0\]
where $\beta$ is skew-symmetric.\\

We will prove\\

\noindent {\bf Theorem 2.}  Locally on $S$, $\xi$ has a special Grothendieck
complex.\\

\noindent {\bf Proof that Theorem 2 $\Rightarrow$ Theorem 1.}\\

We have  \[\psi(\xi_s)=\sum_{i \mbox{ {\scriptsize even}}}\;\mbox{rank }K^i -
\mbox{ rank }(\beta(s))-\sum_{i=0}^m \mbox{ rank }(\alpha^i(s))+ \mbox{ rank
}(\check a^i(s)).\] As $\beta$ is skew-symmetric the parity of $\psi(\xi_s)$
the
same as that of $\displaystyle{\sum_{i\mbox{ {\scriptsize even}}}}$ rank $K^i$
which is (locally) constant.\\

\noindent {\bf \S2.  Special Complexes.}\\

Let $L^*$ be a complex.\\

\noindent Then $L^* \otimes L^*$ is the complex \[(L^* \otimes L^*)^n =
\bigoplus_{n_1 + n_2 = n} L^{n_{1}} \otimes L^{n_{2}}\] with differential
\[d(a\otimes b) = da \otimes b+(-1)^i a \otimes db\] if $a\in L^i$ and $b\in
L^j$.\\

This complex has an involution $\tau:L^*\rightarrow L^*$ given by
$\tau(a\otimes
b) = (-1)^{\alpha \beta} b\otimes a$ where $a\in L^\alpha$ and $b\in
L^\beta$.\\

Let $L^*$ be a complex of free coherent sheaves on $S$.  We will assume
\[L^*:0\rightarrow L^0\rightarrow \;\;\;\;\;\rightarrow L^n\rightarrow 0.\]
Assume that we have a pairing \[\gamma:L^*\otimes L^*\rightarrow {\cal
O}_S(-n)\]
such that $\gamma \circ \tau = (-1)^m\gamma$ and such that $\gamma (s):L^*|_s
\otimes L^*|_s\rightarrow k(s)(-n)$ defines an isomorphism \[L^*|_s\rightarrow
\mbox{ Hom}(L^*|_s, k(s)(-n)).\]

\noindent {\bf Lemma 3.}  $L^*$ is a special complex in a neighborhood of
$s$.\\

\noindent {\it Proof.}  Let $\gamma(a\otimes b)=R(a\otimes b)\cdot 1(-n)$\\

\noindent where $a\in L^p$ and $b \in L^{n-p}$.\\

\noindent Then $R:L^p \otimes L^{n-p} \rightarrow  {\cal O}_S$ satisfies
\[R(\delta_p(\alpha) \otimes \beta)+(-1)^p R(\alpha \otimes
\delta_{n-p-1}(\beta))=0\] for $\alpha \in L^p$ and $\beta \in L^{n-p+1}$.\\

Thus we have a commutative diagram \[\begin{array}{ccc}
L^p&\stackrel{\alpha_p}{\rightarrow }&L^{p+1}\bigskip\\  \downarrow \bar
R_p&&\downarrow \bar R_{p+1}\bigskip\\  \check L^{n-p}&\rightarrow &\check
L^{n-p-1}\\ &(-1)^{p+1}\check \alpha_{n-p-1}.\end{array}\]

If $\bar R_*$ are isomorphisms at $s$, they are isomorphisms in a neighborhood
of $s$.  Hence the high differentials in $L^*$ are isomorphic to the dual of
the
low differential upto sign.\\

We want to check that $R_{m+1}\circ \alpha_m:L^m \rightarrow \check L^m$ is
skew-symmetric.  This will follow if $\bar R_m = \check {\bar R}_{m+1}(-1)^m$.
Now $\check {\bar R}_{m+1}(k)(c) = R(k\otimes c)$ and $\bar R_m(c)(k) =
R(c\otimes k)$.\\

Thus our symmetry condition implies that $\bar R_m=(-1)^{m(m+1)+m}\check
R_{m+1}(= (-1)^m \check R_{m+1})$.  Q.E.D.\\

\noindent {\bf \S4.  The first step.}\\

We fix a point $s$ of $S$ and freely replace $S$ be an open neighborhood of
$s$.  So we may assume that $S$ is affine.  We have the \v Cech resolution
\[\xi
\rightarrow \check {\cal C}^*\] of $\xi$ with respect to some open affine cover
of $X$.  Then $K^*=f_*\check{\cal C}^*$ has homology sheaves $R^if_*\xi$.\\

Now we have a resolution of $L\otimes L$ and $a$ commutative diagram
\[\begin{array}{ccc} \xi \otimes \xi&\rightarrow &\check{\cal C}^*\otimes
\check{\cal C}^*\bigskip \\ \downarrow B&&\downarrow B^*\bigskip \\
\omega_{X/S}&\rightarrow &I^*\end{array}\] where $I^*:0\rightarrow
I^0\rightarrow \cdots\rightarrow I^{2n+1}\rightarrow 0$ where $I^i$ is
injective
if $i<2n+1$ which is a resolution of $\omega_{X/S}$.\\

So we have an induced mapping \[\alpha:K^*\otimes K^* \rightarrow  f_*I^*.\]

We want to replace $\alpha:K^* \otimes K^*\rightarrow f_*I^*$ by a finite
approximation.\\

By Grothendieck's approximation theorem we can find a complex $0\rightarrow
L^0\rightarrow \cdots \rightarrow L^n\rightarrow 0$ of free ${\cal O}_S$-module
of finite type together with a homomorphism $\sigma:L^*\rightarrow K^*$ such
$\sigma$ is a quasi-isomorphism.\\

Thus we get $\beta = \alpha (\sigma\otimes \sigma):L^*\otimes L^* \rightarrow
f_*I^*$.\\

Using Grothendieck's proof we can find a complex $M^*$ of the same kind such
that we have a commutative diagram \[\begin{array}{ccc} L^*\otimes
L^*&\stackrel{i}{\rightarrow }&M^*\\ \searrow&&\downarrow \rho\\
\beta&&f_*I\end{array}\] where $\rho$ is a quasi-isomorphism.\\

Let $i' = (i+(-1)^m i(\tau))/2$.\\

We need another kind of approximation.\\

\noindent {\bf \S5.  The second step.}\\

Let $0\rightarrow L^0\rightarrow \cdots \rightarrow L^p\rightarrow 0$ be a
complex of free coherent sheaves on $S$.  Then replace $S$ by a neighbor of $s$
we may find normalized such complexes $M^*$ and $N^*$ together with
quasi-isomorphisms $M^*\rightarrow L^*\rightarrow N^*$ where normalized means
that the differential of the complex vanishes at $s$.\\

One we do this we'll do the following.\\

We can consider the composition \[R^*\otimes R^* \rightarrow  L^* \otimes L^*
\rightarrow  M^* \rightarrow N^* \stackrel{\alpha}{\rightarrow } {\cal
O}_S(-n)\]
where $R^*$ and $N^*$ are normalized and $\alpha$ is the isomorphism
$n$-homology of $N^*$ same of $M\approx$ same of $f_*i\approx
R^nf_*(\omega_{X/S})\approx {\cal O}_S$ which works on $N^p =0$ if $p>n$ by
similar reasoning.\\

The proof of this step is easy.\\

\noindent Let $f_{i,i}, \ldots, f_{i,i}$ be elements of $L^i$ such that their
reduce in $L^i(s)$ are a basis of a maximal space which is mapped
isomorphically
into $L^{i+1}(s)$.\\

\noindent Then we have a complex $S^*$ such that $S^i$ has basis
$\overline{f_{*,i}}$ and $\overline{df_{*,i-1}}$ where $a_i(f_{*,i} = \overline
{df_{*,i}}$ and $\alpha_{c+1}(df_{*,i})=0$.\\

\noindent Thus we have an obvious homomorphism $S^* \rightarrow L^*$.\\

Let $M^*=L^*/S^*$ is so far quotient.\\

\noindent For $M^*=\check N^*$ where $\check L^*\rightarrow N^*$ is constructed
as before.\\

\noindent {\bf \S6.  The third step.}\\

Now take a $h:R^*\rightarrow K^*$ be a quasi-isomorphism where $R^*$ is
normalized and $k:M^*\rightarrow S^*$ be a quasi-isomorphism where $S^*$ is
normalized.\\

\noindent Then $j'=k\circ i' \circ h:R^* \otimes R^* \rightarrow S^*$.\\

\noindent Let $m:S^* \rightarrow {\cal O}_S(-n)$ be projection on the $n$-the
component.\\

\noindent Then we have $m\circ j':R^* \otimes R^* \rightarrow {\cal
O}_S(-n)$.\\

We want to check that this pairing satisfies the condition of section 2.\\

Clearly $m\circ j' (\tau) = (-1)^m m \otimes j'$ by construction of $i'$.  We
need to check that  \[u^i:R^i(s) \otimes R^{n-i} (s) \rightarrow  k\] is a
perfect pairing.\\

Now by construction $u^i$ is isomorphic $P^i:H^i(X_s, \xi_s) \otimes H^{n-1}
(X_s, \xi_s)\rightarrow k$.\\

To check that this is a perfect pairing by Serre duality it will be enough to
check that it is the usual mapping induced by $m\otimes k \circ i \circ
h|_s$.\\

This is just that \[\rho(\alpha \otimes \beta )=(-1)^{i(n-i)+m} (\rho (\beta
\otimes \alpha))\] but this follows from the assumption on the symmetry of the
pairing $\xi \otimes \xi \rightarrow  \Omega_{X/S}$.\\

This finishes the proof.

\end{document}